\documentclass[a4paper]{article}

\usepackage{url}
\usepackage{longtable}
\usepackage{philex}
\usepackage{amssymb}
\usepackage{latexsym}
\usepackage{bussproofs}
\usepackage{pdflscape}
\usepackage{mathabx}
\usepackage{stackengine}
\stackMath

\usepackage{natbib}

\newtheorem{lemma}{Lemma}
\newtheorem{theorem}{Theorem}

\usepackage{pxfonts}
\usepackage{microtype}

\title{A Binary Quantifier for Definite Descriptions for Cut Free Free Logics}
\author{Nils K\"urbis}
\date{}

\begin{document}
\maketitle

\begin{center}
Published in \emph{Studia Logica}\\
\url{http://dx.doi.org/10.1007/s11225-021-09958-x}\bigskip
\end{center}

\begin{abstract}
\noindent This paper presents rules in sequent calculus for a binary quantifier $I$ to formalise definite descriptions: $Ix[F, G]$ means `The $F$ is $G$'. The rules are suitable to be added to a system of positive free logic. The paper extends the proof of a cut elimination theorem for this system by Indrzejczak by proving the cases for the rules of $I$. There are also brief comparisons of the present approach to the more common one that formalises definite descriptions with a term forming operator. In the final section rules for $I$ for negative free and classical logic are also mentioned. 
\end{abstract}

\noindent Keywords: definite descriptions, free logic, sequent calculus, cut elimination 

\section{Introduction} 
Russell's analysis of the definite article `the' and the ensuing theory of definite descriptions is celebrated as a paradigm for philosophy.\footnote{The phrase stems from \citet[1n]{ramseyphilosophy}. I would like to thank Andrzej Indrzejczak for his comments on this paper and a referee for \emph{Studia Logica}, who also made helpful suggestions for improvement. This paper was written while I was an Alexander von Humboldt Fellow at the University of Bochum. To both institutions many thanks are due.} Consequently, definite descriptions have been the subject of extensive logical and philosophical investigations. It is the more surprising that the formalisation of the theory of definite descriptions has received comparatively little attention when it comes to another paradigm of logic, that of Gentzen's sequent calculus and natural deduction in proof theory. It is almost exclusively due to Andrzej Indrzejczak that the task of combing the two paradigms has been taken on at all. Indrzejczak has provided formalisations of various theories of definite descriptions within sequent calculus, modal as well as non-modal, and proved cut elimination theorems for them \citep{andrzejmodaldescription, andrzejfregean, andrzejexistencedefinedness, andrzejfreedefdescr}. Earlier work in this framework was done by \cite{czermakdefdescr}, and more recently \cite{gratzldefdescr} has formalised Russell's theory of definite descriptions in sequent calculus. 

Most theories of definite descriptions follow Russell in formalising them by means of a term forming operator: $\iota$ binds a variable and forms a singular term out of an open formula. $\iota xFx$ means `the $F$'. The latter are expressions in the grammatical category of names of objects and used accordingly. $G(\iota xFx)$ means `The $F$ is $G$'. Few axiomatisations of theories of definite descriptions, however, follow Russell in some other respects. Russell considered definite descriptions to be incomplete symbols. The meaning of $\iota xFx$ is given by a contextual definition and it has no meaning outside the context of a formula in which it occurs. \citep[69ff]{russellwhitehead} Indeed, upon analysis, the definite description disappears altogether: `The $F$ is $G$' means no more nor less than $\exists y(\forall x(Fx\leftrightarrow x=y)\land Gy)$. For Russell, the use of the $\iota$ operator is a matter of convenience, as it can be dispensed with. A formula with the $\iota$ operator is an abbreviation of a longer formula and as such often easier to read. But that convenience is eradicated again by the need Russell saw for marking scope distinctions: $\neg G(\iota xFx)$ is ambiguous, as it may represent either the internal negation of `The $F$ is $G$', that is `The $F$ is not $G$', or its external negation, `It is not the case that the $F$ is $G$'. Russell avoids ambiguity with a rather clumsy method of marking scope, which consists in repeating the $\iota$ term in square brackets at the beginning of its scope. The internal negation of `The $F$ is $G$' is formalised as $[\iota xF]\neg G(\iota xFx)$, its external negation as $\neg [\iota xF]G(\iota xF)$. \citep[181ff]{russellwhitehead} It is fair to say that these two formulas lose much of the greater ease of readability that $G(\iota xFx)$ might have over $\exists y(\forall x(F\leftrightarrow x=y)\land Gy)$. 

The most common formalisations of theories of definite descriptions do not follow Russell in incorporating scope distinctions. In the classic work of Lambert, $\iota$ is axiomatised on the basis of what is now known as \emph{Lambert's Law}:\footnote{A principle almost like Lambert's Law was proposed by \cite{hintikkatowardsdefdesc}, but subsequently shown to be inconsistent by \cite{lambertnotesEIII}: Hintikka omitted the outermost quantifier. The latter paper also contains the first formulation by Lambert of Lambert's Law.}

\lbp{LL}{$LL$}{$\forall y(\iota xFx=y \leftrightarrow \forall x(Fx\leftrightarrow x=y))$}

\noindent which makes no mention of scope distinctions. \rf{LL} axiomatises what is commonly regarded as the minimal theory of $\iota$. The logic can be either a negative or a positive free logic, with many free logicians preferring the latter. Additional axioms for stronger theories considered by \cite{lambertnotesEIV}, \cite{fraassenthexxlambert}, \cite{bencivengahandbook} and others also do not provide means for distinguishing scope.\footnote{For an overview see, besides the articles just quoted, \cite{morschersimonsfreelogic}.} 

It is the expressed aim of formalisations of theories of definite descriptions following Hintikka and Lambert that the theory should only care for the proper definite descriptions, that is, the cases where there is a unique object that satisfies the predicate $F$ in $\iota xFx$, and remain largely silent if there is not. Hintikka makes the point that `there is little to be said about the properties of a described object unless we know that it exists', but whether it exists or not, to be the $F$ means to be a unique $F$ \citep[83]{hintikkatowardsdefdesc}. Hintikka's theory still said rather too much about definite descriptions -- it was inconsistent -- but, \emph{mutatis mutandis}, his observation motivates Lambert's Law. Lambert concurs in spirit \citep[2f]{lambertnotesEII}. Van Fraassen, too, underlines the neutrality  of the minimal theory of definite descriptions when it comes to improper definite descriptions \citep[9f]{fraassenthexxlambert}. Bencivenga notes that the motivation behind it is that `everybody agrees on how to treat denoting descriptions, and [the minimal theory of definite descriptions] says nothing (specific) beyond that' \citep[417]{bencivengahandbook}. But if there is a unique $F$, then, as is also the case in Russell's theory, scope distinctions no longer matter, and so there is no need for them in this theory. 

Exceptions to the rule are provided by Lambert himself in his formalisation of a Russellian theory of definite descriptions within negative free logic, where scope is marked by an operator for predicate abstraction \citep{lambertfreedef}. This method is rather more elegant than Russell's own. It is also used by \cite{mendelsohnfitting} and \cite{garsonmodallogic} in their investigations of definite descriptions in modal extensions of positive free logic. There is thus a place for a means for marking scope distinctions also in positive free logic.

An intriguing alternative is to formalise sentences containing definite descriptions by a binary quantifier which incorporates scope distinctions directly into the notation. This approach was recommended by Dummett , who proposes that `The $F$ is $G$' should be formalised by an expression $Ix[F, G]$, where $I$ binds a variable and forms a formula out of two formulas \citep[p.162]{dummettfregelanguage}.\footnote{A closely related notation is used by \cite{nealedescriptions} and briefly by \citet[Sec. 8.4]{bostockintermediate}.} The proof theory of $I$ was investigated within natural deduction for intuitionist negative and positive free logics in \citep{kurbisiotaI, kurbisiotaII, kurbisiotaIII}. The present paper investigates the proof theory of $I$ in sequent calculus for classical positive free logic. I will give rules for $I$ suitable to this framework and briefly compare them to axioms for $\iota$. The main part of the paper consists in a proof of a cut elimination theorem for the resulting system. It builds on a result of Indrzejczak's published recently in this journal \citep{andrzejcutfreefreelogic}. Indrzejczak proves cut elimination theorems for a variety of positive and negative free logics. In the present paper Indrzejczak's proof is continued by the cases covering the rules for $I$. The contribution of this paper is thus two-fold: to propose a formalisation of a theory of definite descriptions within classical positive free logic that incorporates a means for marking scope distinctions, and to show that this is done in a proof-theoretically satisfactory way. 

The generality of Indrzejczak's result means that one could envisage adding the rules for $I$ studied here to other systems of free logic. However, different rules may be better suited to different logics. In particular, in negative free logic significantly simpler rules for $I$ will do than those suitable for positive free logic proposed here. In the final section I consider them briefly. Indrzejczak's cut elimination theorem could be extended to cover negative free logic extended by $I$. The resulting system provides a proof-theoretically satisfactory formalisation of a Russellian theory of definite descriptions.

\section{A System of Positive Free Logic} 
The language is standard. Free variables are distinguished from bound ones by the use of parameters $a, b, c\ldots$ for the former and $x, y, z \ldots$ for the latter. For the purposes of the present section, the terms of the language are the parameters, constants and complex terms formed from them by function symbols. The latter play virtually no role in the present paper, except briefly in the conclusion, but as complex terms formed by the $\iota$ operator are of course prominent in the other sections, we might as well include function symbols here. The occurrence of free variables in formulas will not be indicated explicitly except where substitution is concerned. Instead of $A(x)$, I'll simply write $A$. $A_t^x$ is the result of substituting $t$ for $x$ in $A$, where it is assumed that no variable free in $t$ becomes bound in $A_t^x$, i.e. that $t$ is free for $x$ in $A$. Exceptions are the existence predicate, where I'll write $\exists !x$ and $\exists !t$, and in the following section, as in the previous one, I'll continue to use brackets  where substitution of variables by $\iota$ terms are concerned, as in $G(\iota xF)$, except, again, where $G$ is the existence predicate. In all cases, identities are written as usual.  

 $\Gamma, \Delta$ denote finite multisets of formulas. Indrzejczak's system GPFL has the following rules:\footnote{It is an extension of the propositional G1 calculus of \citet[61f]{troelstraschwichtenberg} by quantifier rules suitable to free logic and standard rules for identity.}
 
\def\fCenter{\ \Rightarrow\ }

\begin{longtable}{l l}
(Ax) \ \AxiomC{$A\fCenter A$}
\DisplayProof & 
\AxiomC{$\Gamma\Rightarrow\Theta, A$}
\AxiomC{$A, \Delta\Rightarrow \Lambda$}
\LeftLabel{Cut \ }
\BinaryInfC{$\Gamma, \Delta\Rightarrow\Theta, \Lambda$}
\DisplayProof\\
 \\
\Axiom$\Gamma\fCenter \Delta$
\LeftLabel{(LW)\ }
\UnaryInf$A, \Gamma\fCenter \Delta$
\DisplayProof  &
\Axiom$\Gamma\fCenter \Delta$
\LeftLabel{(RW) \ }
\UnaryInf$\Gamma\fCenter \Delta, A$
\DisplayProof\\
\\
\Axiom$A, A, \Gamma\fCenter \Delta$
\LeftLabel{(LC) \ }
\UnaryInf$A, \Gamma\fCenter \Delta$
\DisplayProof & 
\Axiom$\Gamma\fCenter \Delta, A, A$
\LeftLabel{(RC) \ }
\UnaryInf$\Gamma\fCenter\Delta, A$
\DisplayProof\\
\\
\Axiom$\Gamma\fCenter \Delta, A$
\LeftLabel{$(L\neg)$ \ }
\UnaryInf$\neg A, \Gamma\fCenter \Delta$
\DisplayProof &
\Axiom$A, \Gamma\fCenter \Delta$
\LeftLabel{$(R\neg)$ \ }
\UnaryInf$\Gamma\fCenter \Delta, \neg A$
\DisplayProof\\
\\
\Axiom$A, B, \Gamma\fCenter \Delta$
\LeftLabel{$(L\land)$ \ }
\UnaryInf$A\land B, \Gamma\fCenter \Delta$
\DisplayProof &
\AxiomC{$\Gamma\Rightarrow \Delta, A$} 
\AxiomC{$\Gamma\Rightarrow \Delta, B$}
\LeftLabel{$(R\land)$ \ }
\BinaryInfC{$\Gamma \Rightarrow \Delta,  A\land B$}
\DisplayProof \\
\\
\AxiomC{$A, \Gamma\Rightarrow \Delta$}
\AxiomC{$B, \Gamma\Rightarrow \Delta$}
\LeftLabel{$(L\lor)$ \ }
\BinaryInfC{$A\lor B, \Gamma \Rightarrow \Delta$}
\DisplayProof &
\Axiom$\Gamma\fCenter\Delta, A, B$
\LeftLabel{$(R\lor)$ \ }
\UnaryInf$\Gamma\fCenter \Delta, A\lor B$
\DisplayProof\\
\\
\AxiomC{$\Gamma\Rightarrow \Delta, A$} 
\AxiomC{$B, \Gamma\Rightarrow \Delta$}
\LeftLabel{$(L\!\rightarrow)$ \ }
\BinaryInfC{$A\rightarrow B, \Gamma\Rightarrow \Delta$}
\DisplayProof & 
\Axiom$A, \Gamma\fCenter \Delta, B$
\LeftLabel{$(R\!\rightarrow)$ \ }
\UnaryInf$\Gamma\fCenter \Delta, A\rightarrow B$
\DisplayProof\\
\\
\Axiom$A_t^x, \Gamma\fCenter \Delta$
\LeftLabel{$(L\forall)$ \ }
\UnaryInf$\exists !t, \forall xA, \Gamma \fCenter \Delta$
\DisplayProof & 
\Axiom$\exists ! a, \Gamma \fCenter \Delta, A_a^x$
\LeftLabel{$(R\forall)$ \ }
\UnaryInf$\Gamma\fCenter \Delta, \forall xA$
\DisplayProof\\
\\
\Axiom$\exists ! a, A_a^x, \Gamma\fCenter \Delta$
\LeftLabel{$(L\exists)$ \ }
\UnaryInf$\exists xA, \Gamma\fCenter \Delta$
\DisplayProof & 
\AxiomC{$\Gamma\Rightarrow\Delta, A_t^x$}
\LeftLabel{$(R\exists)$ \ }
\UnaryInfC{$\exists ! t,\Gamma \Rightarrow \Delta, \exists xA$}
\DisplayProof
\end{longtable}

\noindent where in $(L\exists)$ and $(R\forall)$, $a$ does not occur in the conclusion.  

Indrzejczak's system GPFL$_=$ is formed by adding rules for identity to GPFL:

\begin{center} 
\Axiom$A_{t_2}^x, \Gamma\fCenter\Delta$
\LeftLabel{$(=I)$ \ }
\UnaryInf$t_1=t_2, A_{t_1}^x, \Gamma\fCenter\Delta$
\DisplayProof\qquad
\Axiom$t=t, \Gamma\fCenter\Delta$
\LeftLabel{$(=E)$ \ }
\UnaryInf$\Gamma\fCenter\Delta$
\DisplayProof
\end{center} 

\noindent where $A$ is atomic. The general case follows by induction. 

Indrzejczak proves that cut is eliminable from GPFL and GPFL$_=$ \citep[Theorem 3]{andrzejcutfreefreelogic}. In the next section I will extend GPFL$_=$ by rules for the binary quantifier $I$ and in the section thereafter continue Indrzejczak's proof to show that cut is eliminable also from the resulting system GPFL$_=^I$. 

For comparisons between GPFL$_=^I$ and a system with the term forming $\iota$ operator, it will be useful to have rules for the biconditional: 

\begin{center} 
\AxiomC{$\Gamma\Rightarrow \Delta, A, B$} 
\AxiomC{$A, B, \Gamma\Rightarrow \Delta$}
\LeftLabel{$(L\leftrightarrow)$ \ }
\BinaryInfC{$A\leftrightarrow B, \Gamma\Rightarrow \Delta$}
\DisplayProof\bigskip

\AxiomC{$A, \Gamma\Rightarrow \Delta, B$}
\AxiomC{$B, \Gamma\Rightarrow \Delta, A$}
\LeftLabel{$(R\leftrightarrow)$ \ }
\BinaryInfC{$\Gamma\Rightarrow \Delta, A\leftrightarrow B$}
\DisplayProof
\end{center} 

\noindent These are derivable from the rules for $\rightarrow$ and $\land$ taking the usual definition of $\leftrightarrow$. 

Two useful provable sequents are $A, A\leftrightarrow B\Rightarrow B$ and $A\leftrightarrow B, B\Rightarrow A$, which I will call $(MP\leftrightarrow)$. The first is proved in the following way, the second similarly: 

\begin{center}
\Axiom$A\fCenter A$
\doubleLine
\UnaryInf$A\fCenter B, A, B$
\Axiom$B\fCenter B$
\doubleLine
\UnaryInf$A, B, A\fCenter B$
\BinaryInf$A \leftrightarrow B, A\fCenter B$
\DisplayProof
\end{center}

\noindent Here and in the following, double lines indicate possibly multiple applications of rules, in particular the structural rule weakening, which must be used abundantly to make the contexts of the operational rules identical.

\section{Adding $I$}
The syntax of $\iota$ is that if $F$ is a formula, $\iota xF$ is a term. The syntax of $I$ is that if $F$ and $G$ are formulas, $Ix[F, G]$ is a formula. In both cases $x$ is bound. 

Let GPFL$_=^\iota$ be GPFL$_=$ with its language extended by $\iota$ and \rf{LL} added as an axiom. In this system, what we might call the primary occurrences of $\iota$ terms are those where they occupy the left or right of $=$. Occurrences where a predicate $G$ other than identity is applied to an $\iota$ term are secondary: the logic of $G(\iota xF)$ is explained in terms of and derived from primary occurrences of $\iota xF$. But it would be possible to start the other way round. The following two principles are jointly equivalent to \rf{LL} in GPFL$_=^\iota$:\footnote{This holds already in intuitionist positive free logic. For proof see \citep[ $\ast 3$, $\ast 4$]{kurbisiotaIII}.}

\lbp{L1'}{$\iota1$}{$\exists y(\forall x (F\leftrightarrow x=y)\land G)\rightarrow G(\iota xF)$}

\lbp{L2'}{$\iota2$}{$(G(\iota xF)\land \exists !\iota xF) \rightarrow \exists y(\forall x (F\leftrightarrow x=y)\land G)$}

\noindent In positive free logic, $\exists y(\forall x (F\leftrightarrow x=y))$ is equivalent to $\exists !\iota x F$, so the Russellian phrase $\exists y(\forall x (F\leftrightarrow x=y)\land G)$ is equivalent to $G(\iota xF)\land \exists !\iota xF$. 

Let GPFL$_=^I$ be GPFL$_=$ with its language extended by $I$ and these rules added:

\bigskip

\noindent
\AxiomC{$\Gamma\Rightarrow \Delta, F_t^x$}
\AxiomC{$\Gamma\Rightarrow \Delta, G_t^x$}
\AxiomC{$\Gamma\Rightarrow \Delta, \exists !t$}
\AxiomC{$\exists ! a, F_a^x, \Gamma \Rightarrow \Delta, a=t$}
\LeftLabel{$(RI)$ \ }
\QuaternaryInfC{$\Gamma \Rightarrow \Delta, I x[F, G]$}
\DisplayProof

\bigskip

\noindent where $a$ does not occur in the conclusion. 

\bigskip

\noindent {\small
\AxiomC{$\Gamma\Rightarrow \Delta, F_t^x$}
\AxiomC{$\Gamma\Rightarrow \Delta, \exists !t$}
\AxiomC{$F_a^x, \exists ! a, \Gamma\Rightarrow \Delta, a=t$}
\AxiomC{$F_b^x, G_b^x, \exists !b, \Gamma\Rightarrow \Delta$}
\LeftLabel{$(LI^1)$ \ }
\QuaternaryInfC{$Ix[F, G], \Gamma\Rightarrow \Delta$}
\DisplayProof
}

\bigskip

\noindent where $a$ and $b$ do not occur in the conclusion. 

\bigskip

\noindent {\small
\AxiomC{$\Gamma\Rightarrow \Delta, F_{t_1}^x$}
\AxiomC{$\Gamma\Rightarrow \Delta, F_{t_2}^x$}
\AxiomC{$\Gamma\Rightarrow \Delta, \exists !t_1$}
\AxiomC{$\Gamma\Rightarrow \Delta, \exists !t_2$}
\AxiomC{$\Gamma\Rightarrow \Delta, A_{t_2}^x$}
\LeftLabel{$(LI^2)$ \ }
\QuinaryInfC{$Ix[F, \exists !x], \Gamma\Rightarrow \Delta, A_{t_1}^x$}
\DisplayProof
}

\bigskip

\noindent where $A$ is an atomic formula.

\bigskip
\noindent
\Axiom$F_a^x, \exists !a, \Gamma\fCenter \Delta$
\LeftLabel{$(LI^3)$ \ }
\UnaryInf$Ix[F, \exists !x], \Gamma\fCenter\Delta$
\DisplayProof

\bigskip

\noindent where $a$ does not occur in the conclusion.

\bigskip

\noindent
\AxiomC{$\Gamma\Rightarrow \Delta, F_{t_1}^x$}
\AxiomC{$\Gamma\Rightarrow \Delta, \exists !t_1$}
\AxiomC{$\Gamma\Rightarrow \Delta, \exists !t_2$}
\AxiomC{$\Gamma\Rightarrow \Delta, A_{t_2}^x$}
\LeftLabel{$(LI^4)$ \ }
\QuaternaryInfC{$Ix[F, x=t_2], \Gamma\Rightarrow \Delta, A_{t_1}^x$}
\DisplayProof

\bigskip

\noindent where $A$ is an atomic formula. 

\bigskip

\noindent
\AxiomC{$\Gamma\Rightarrow \Delta, \exists !t$}
\AxiomC{$F_a^x, \exists !a, \Gamma\Rightarrow \Delta$}
\LeftLabel{$(LI^5)$ \ }
\BinaryInfC{$Ix[F, x=t], \Gamma\Rightarrow\Delta$}
\DisplayProof

\bigskip

\noindent where $a$ does not occur in the conclusion. 

\bigskip

\noindent These rules are those of \citep{kurbisiotaIII} transposed to sequent calculus. That paper also contains an extensive discussion of the conceptual foundations of the present formalisation of definite descriptions and explains why these rules are adequate for the account at hand. Here I only note two things. First, $(LI^4)$ and $(LI^5)$ are required to mimic some inferences in the framework using $\iota$ involving identity, of which I will give an example shortly. Secondly, the remaining rules are equivalent to principles corresponding to \rf{L1'} and \rf{L2'} under a translation between the languages of GPFL$^\iota_=$ and GPFL$^I_=$ in which $G(\iota xF)$, $\exists !\iota xF$, $\iota xF=t$ are translated as $Ix[F, G]$, $Ix[F, \exists !x]$, $Ix[F, x=t]$, respectively:\footnote{A more precise account of this translation may be found in \citep{kurbisiotaII}.} 

\lbp{L1}{$SI1$}{$\exists y(\forall x (F\leftrightarrow x=y)\land G_y^x)\Rightarrow Ix[F, G]$}

\lbp{L2}{$SI2$}{$Ix [F, G], Ix[F, \exists !x]\Rightarrow \exists y(\forall x (F\leftrightarrow x=y)\land G_y^x)$}

\noindent For simplicity I will use a more convenient, but equivalent, version of $(LI^2)$: 

\bigskip

\noindent
\AxiomC{$\Gamma\Rightarrow \Delta, F_{t_1}^x$}
\AxiomC{$\Gamma\Rightarrow \Delta, F_{t_2}^x$}
\AxiomC{$\Gamma\Rightarrow \Delta, \exists !t_1$}
\AxiomC{$\Gamma\Rightarrow \Delta, \exists !t_2$}
\LeftLabel{$(LI{^2}')$ \ }
\QuaternaryInfC{$Ix[F, \exists !x], \Gamma\Rightarrow \Delta, t_1=t_2$}
\DisplayProof

\bigskip

\noindent The reason $(LI{^2}')$ is not the official rule of the present system is that with it, cuts on identities that are concluded by $(LI{^2}')$ in the left premises and by $(=E)$ in the right premise are not eliminable.\footnote{\label{andrzejfootnote}Indrzejczak suggests in correspondence that this problem can be avoided with an alternative to $(=I)$: from $\Gamma\Rightarrow \Delta, t_1=t_2$ and $\Gamma\Rightarrow \Delta, A_{t_1}$ infer $\Gamma\Rightarrow \Delta, A_{t_2}$. He also proposes a further version of $(LI^2)$ which avoids the problem while keeping the original rule $(=I)$: from $\Gamma\Rightarrow \Delta, F_{t_1}^x$, $\Gamma\Rightarrow \Delta, F_{t_2}^x$, $\Gamma\Rightarrow \Delta, \exists !t_1$, $\Gamma\Rightarrow \Delta, \exists !t_2$ and $t_1=t_2, \Gamma \Rightarrow \Delta$ infer $Ix[F, \exists !x], \Gamma\Rightarrow \Delta$. Similarly for $(LI^4)$. My aim here is to stay close to the system as presented in his paper, and I mention the first option only for its interest. To the second one I'll come back in Section 5.1.}

\begin{theorem}\label{equiv}
\normalfont
Given the rules of GPFL$_=$: 

\noindent (a) \rf{L1} and $(RI)$ are interderivable; 

\noindent (b) $(LI^1)$ is derivable from \rf{L2} and the instance of \rf{L1} with $G$ replaced by $\exists !$; 

\noindent (c) $(LI{^2}')$ and $(LI^3)$ are derivable from the instance of \rf{L2} with $G$ replaced by $\exists!$; 

\noindent (d) \rf{L2} is derivable from $(LI^1)$, $(LI^2)$ and $(LI^3)$. 
\end{theorem}

\noindent \emph{Proof.}\bigskip

\noindent (a.i) Assume sequents (1) $\Gamma\Rightarrow \Delta, F_t^x$, (2) $\Gamma\Rightarrow \Delta, G_t^x$, (3) $\Gamma\Rightarrow \Delta, \exists !t$ and (4) $\exists ! a, F_a^x, \Gamma \Rightarrow \Delta, a=t$, where $a$ is not free in $\Gamma, \Delta$. From $F_a^x\Rightarrow F_a^x$ by $(=I)$ we have $a=t, F_t^x\Rightarrow F_a^x$, and so from (1) by Cut $a=t, \Gamma\Rightarrow \Delta, F_a^x$. Then from (4) by weakening and $(R\leftrightarrow)$ $\exists ! a, \Gamma \Rightarrow \Delta, F_a^x\leftrightarrow a=t$, and so by $(R\forall)$ we derive $\Gamma\Rightarrow \Delta, \forall x(F\leftrightarrow x=t)$. So from (2) by $(R\land)$: $\Gamma \Rightarrow \Delta, \forall x(F\leftrightarrow x=t)\land G_t^x$, and by $(R\exists)$: $\exists ! t, \Gamma\Rightarrow \Delta, \exists y(\forall x(F\leftrightarrow x=y)\land G_y^x)$, so from (3) by Cut and contraction, $\Gamma\Rightarrow \Delta, \exists y(\forall x(F\leftrightarrow x=y)\land G_y^x)$. Finally, by Cut from \rf{L1} we conclude $\Gamma\Rightarrow \Delta, Ix[F, G]$. 

\bigskip 

\noindent (a.ii) First, we prove two sequents using $(MP\leftrightarrow)$, $(L\forall)$ and $(=E)$:\bigskip

\AxiomC{$b=b, F_b^x\leftrightarrow b=b \Rightarrow F_b^x$}
\UnaryInfC{$b=b, \exists !b, \forall x(F\leftrightarrow x=b)\Rightarrow F_b^x$}
\UnaryInfC{$\exists !b, \forall x(F\leftrightarrow x=b)\Rightarrow F_b^x$}
\DisplayProof\qquad
\AxiomC{$F_a^x, F_a^x\leftrightarrow a=b \Rightarrow a=b$}
\UnaryInfC{$\exists ! a, F_a^x, \forall x(F\leftrightarrow x=b)\Rightarrow a=b$}
\DisplayProof

\bigskip

\noindent Then from these and two axioms by weakening we derive the premises of $(RI)$:\bigskip 

{\footnotesize
\def\defaultHypSeparation{\hskip .05in}
\Axiom$\exists !b, \forall x(F\leftrightarrow x=b)\fCenter F_b^x$
\Axiom$ G_b^x\fCenter G_b^x$
\Axiom$\exists !b\fCenter \exists ! b$
\Axiom$\exists ! a, F_a^x, \forall x(F\leftrightarrow x=b)\fCenter a=b$
\doubleLine
\QuaternaryInf$ \exists !b, \forall x(F\leftrightarrow x=b), G_b^x\fCenter Ix[F, G]$
\UnaryInf$ \exists !b, \forall x(F\leftrightarrow x=b)\land G_b^x\fCenter Ix[F, G]$
\UnaryInf$ \exists y(\forall x(F\leftrightarrow x=y)\land G_y^x)\fCenter Ix[F, G]$
\DisplayProof
}
\bigskip

\noindent (b) Assume (1) $\Gamma\Rightarrow \Delta, F_t^x$, (2) $\Gamma\Rightarrow \Delta, \exists !t$, (3) $F_a^x, \exists ! a, \Gamma\Rightarrow \Delta, a=t$ and (4) $F_b^x, G_b^x, \exists !b, \Gamma\Rightarrow \Delta$, $a$ and $b$ not free in $\Gamma, \Delta$. Using (1), (2) twice and (3), by the derivability of $(RI)$ from \rf{L1} and replacing $G$ with $\exists !$, infer $\Gamma\Rightarrow \Delta, Ix[F, \exists ! x]$. So from \rf{L2} by Cut $Ix [F, G], \Gamma\Rightarrow \Delta, \exists y(\forall x (F\leftrightarrow x=y)\land G_y^x)$. By the rules for $\exists$ and $\land$, $\exists y(\forall x (F\leftrightarrow x=y)\land G_y^x)\Rightarrow \exists x(F\land G)$. From (4) by the same rules $\exists x(F\land G), \Gamma\Rightarrow \Delta$, so by Cut twice and contraction $Ix[F, G], \Gamma\Rightarrow \Delta$. 

\bigskip

\noindent (c) This is fairly straightforward, so it is left as an exercise. 

\bigskip

\noindent (d) This is not so straightforward. Let $\Pi$ be the following deduction, which ends with an application of $(LI{^2}')$ to sequents derived from axioms by weakening: \bigskip

\AxiomC{$F_a^x\Rightarrow F_a^x$}
\AxiomC{$F_b^x\Rightarrow F_b^x$}
\AxiomC{$\exists !a\Rightarrow \exists !a$}
\AxiomC{$\exists !b\Rightarrow \exists !b$}
\doubleLine
\QuaternaryInfC{$Ix[F, \exists ! x], F_a^x, F_b^x, \exists !a, \exists !b\Rightarrow a=b$}
\DisplayProof

\bigskip

\noindent Let $\Sigma$ the following deduction: \bigskip

\AxiomC{$\Pi$}
\Axiom$F_a^x\fCenter F_a^x$
\UnaryInf$a=b, F_b^x\fCenter F_a^x$
\doubleLine
\BinaryInf$Ix[F, \exists ! x], F_b^x, \exists !a, \exists !b\fCenter F_a^x\leftrightarrow a=b$ 
\UnaryInf$Ix[F, \exists !x], F_b^x, \exists ! b\fCenter \forall x(F\leftrightarrow x=b)$
\AxiomC{$G_b^x\Rightarrow G_b^x$}
\doubleLine
\BinaryInfC{$G_b^x, Ix[F, \exists !x], F_b^x, \exists !b\Rightarrow \forall x(F\leftrightarrow x=b)\land G_b^x$} 
\UnaryInfC{$\exists ! b, G_b^x, Ix[F, \exists !x], F_b^x, \exists !b\Rightarrow \exists y(\forall x(F\leftrightarrow x=y)\land G_y^x)$} 
\UnaryInfC{$G_b^x, Ix[F, \exists !x], F_b^x, \exists !b\Rightarrow \exists y(\forall x(F\leftrightarrow x=y)\land G_y^x)$} 
\DisplayProof

\bigskip 

\noindent We now put $\Pi$ with $b$ replaced by a fresh parameter $c$ and $\Sigma$ together with two more premises for an application of $(LI^1)$ after some steps by weakening and continue with $(LI^3)$ and contraction:\bigskip

\Axiom$F_c^x\fCenter F_c^x$
\Axiom$\exists !c\fCenter \exists ! c$
\AxiomC{$\Pi_c^b$}
\AxiomC{$\Sigma$}
\doubleLine
\QuaternaryInf$Ix[F, G], F_c^x, \exists !c,  Ix[F, \exists !x]\fCenter \exists y(\forall x (F\leftrightarrow x=y)\land G_y^x)$
\UnaryInf$Ix[F, G], Ix[F, \exists !x], Ix[F, \exists !x]\fCenter \exists y(\forall x (F\leftrightarrow x=y)\land G_y^x)$
\UnaryInf$Ix[F, G], Ix[F, \exists !x]\fCenter \exists y(\forall x (F\leftrightarrow x=y)\land G_y^x)$
\DisplayProof

\bigskip 

\noindent This completes the proof of theorem \ref{equiv}. 

\bigskip

\noindent To close this section here is a sketch of a proof of an important principle in which $(LI^4)$ and $(LI^5)$ are indispensable.\footnote{This addition was requested by a referee to make the discussion self-contained.} All parameters are fresh, applications of structural rules left implicit. First, the sequent (1) $Ix[F, x=a], \exists !a\Rightarrow F_a^x$ is provable: replacing $t_1$ by $b$, $t_2$ by $a$ and $A_{t_2}$ by $a=a$ establishes the sequent $Ix[F, x=a],  F_b^x, \exists !b, \exists !a \Rightarrow b=a$ by $(LI^4)$ and $(=E)$, whence by $(=I)$ and Cut, $Ix[F, x=a],  F_b^x, \exists !b, \exists !a \Rightarrow F_a^x$. (1) follows by $(LI^5)$. Using once more $(LI^4)$, this time replacing $t_1$ by $c$ and $t_2$ by $t$, proves the sequent (2) $Ix[F, x=t],  F_c^x, \exists !c, \exists !t \Rightarrow c=t$. Using $\exists ! a\Rightarrow \exists ! a$ both as the second and third premises of $(RI)$, (1) as the first, (2) as the fourth, establishes $Ix[F, x=a], \exists !a\Rightarrow Ix[F, \exists !x]$. Finally, an application of $(L\exists)$ derives $\exists yIx[F, x=y]\Rightarrow Ix[F, \exists !x]$. The converse is left as an exercise. Thus `Something is identical to the $F$' is equivalent to `The $F$ exists'. This principle is an aspect where free definite description theorists agree with Russell. It shows that identity sometimes carries aspects of existence, and $(LI^4)$ and $(LI^5)$ ensure that this is also the case in the present formalisation of definite descriptions.

\section{Cut Elimination} 
We continue Indrzejczak's proof of Cut elimination for GPFL$_=$ and check that the Right and Left Reduction Lemmas hold for GPFL$_=^I$ by checking the rules for $I$: consequently Cut elimination holds for the latter system. $d(A)$ is the degree of the formula $A$, that is the number of connectives occurring in it. $\exists ! t$ is atomic, that is of degree $0$. For a proof $\mathcal{D}$, $d(\mathcal{D})$ is the degree of the highest degree of any cut formula in $\mathcal{D}$. $A^k$ denotes $k$ occurrences of $A$, $\Gamma^k$ $k$ occurrences of the formulas in $\Gamma$. The height of a deduction is the largest number of rules applied above the conclusion, that is the number of nodes of the longest branch in the deduction. $\vdash_k \Gamma\Rightarrow \Lambda$ means that the sequent has a proof of height $k$. This is used only in the Substitution Lemma: 

\begin{lemma}
If $\vdash_k \Gamma\Rightarrow \Delta$, then $\vdash_k \Gamma_t^a\Rightarrow \Delta_t^a$. 
\end{lemma} 

\noindent The proof goes through as usual. 

\begin{lemma}[Right Reduction]
If $\mathcal{D}_1\vdash \Theta\Rightarrow \Lambda, A$, where $A$ is principal, and $\mathcal{D}_2\vdash A^k, \Gamma\Rightarrow\Delta$ have degrees $d(\mathcal{D}_1), d(\mathcal{D}_2)<d(A)$, then there is a proof $\mathcal{D}\vdash \Theta^k, \Gamma\Rightarrow \Lambda^k, \Delta$ with $d(\mathcal{D})<d(A)$. 
\end{lemma}

\noindent \emph{Proof.} By induction over the height of $\mathcal{D}_2$. 

The basis is trivial: if $d(\mathcal{D}_2)=1$, then $A^k, \Gamma\Rightarrow\Delta$ is an axiom and hence $k=1$, $\Gamma$ is empty, and $\Delta$ consists of only one $A$, and we need to show $\Theta\Rightarrow \Lambda, A$, but that is already proved by $\mathcal{D}_1$. 

For the induction step, we consider the rules for $I$:\bigskip

\noindent (I) The last step of $\mathcal{D}_2$ is by $(RI)$. Then the occurrences $A^k$ in the conclusion of $\mathcal{D}_2$ are parametric and occur in all four premises of $(RI)$: apply the induction hypothesis to them and apply $(RI)$ afterwards. The result is the desired proof $\mathcal{D}$.\bigskip

\noindent (II) The last step of $\mathcal{D}_2$ is by $(LI^1)$. There are two cases: 

\noindent (II.a) The principal formula $Ix[F, G]$ of $(LI^1)$ is not one of the $A^k$: apply the induction hypothesis to the premises of $(LI^1)$ and then apply the rule. 

\noindent (II.b) The principal formula $Ix[F, G]$ of $(LI^1)$ is one of the $A^k$. Let $\Xi$ be $Ix[F, G]^{k-1}$, i.e. $\Xi$ consists of $k-1$ occurrences of $Ix[F, G]$, then $\mathcal{D}_2$ ends with: 

\bigskip

\begin{center}
\def\defaultHypSeparation{\hskip .05in}
\AxiomC{$\Xi, \Gamma\Rightarrow \Delta, F_t^x$}
\AxiomC{$\Xi, \Gamma\Rightarrow \Delta, \exists !t$}
\AxiomC{$F_a^x, \exists ! a, \Xi, \Gamma\Rightarrow \Delta, a=t$}
\AxiomC{$F_b^x, G_b^x, \exists !b, \Xi, \Gamma\Rightarrow \Delta$}
\QuaternaryInfC{$Ix[F, G]^k, \Gamma\Rightarrow \Delta$}
\DisplayProof
\end{center}

\bigskip

\noindent By induction hypothesis we have: \bigskip

(1) $\Theta^{k-1}, \Gamma\Rightarrow \Lambda^{k-1}, \Delta, F_t^x$ 

(2) $\Theta^{k-1}, \Gamma\Rightarrow \Lambda^{k-1}, \Delta, \exists !t$ 

(3) $F_a^x, \exists ! a, \Theta^{k-1}, \Gamma\Rightarrow \Lambda^{k-1}, \Delta, a=t$ 

(4) $ F_b^x, G_b^x, \exists !b, \Theta^{k-1},\Gamma\Rightarrow\Lambda^{k-1}, \Delta$\bigskip 

\noindent We only need (4), from which by the Substitution Lemma we get:\bigskip 

(5) $\Theta^{k-1}, F_t^x, G_t^x, \exists !t, \Gamma\Rightarrow\Lambda^{k-1}, \Delta$\bigskip

\noindent $A$ is principal in $\mathcal{D}_1$, so it ends with: 

\bigskip

\begin{center} 
\AxiomC{$\Theta\Rightarrow \Lambda, F_t^x$}
\AxiomC{$\Theta\Rightarrow \Lambda, G_t^x$}
\AxiomC{$\Theta\Rightarrow \Lambda, \exists !t$}
\AxiomC{$\exists ! a, F_a^x, \Theta \Rightarrow \Lambda, a=t$}
\QuaternaryInfC{$\Theta \Rightarrow \Lambda, Ix[F, G]$}
\DisplayProof
\end{center} 

\bigskip

\noindent Apply cut three times, to (5) and each of the first three premises, conclude $\Theta^k, \Gamma\Rightarrow\Lambda^k, \Delta$ by contraction.\bigskip 

\noindent (III) The last step of $\mathcal{D}_2$ is by $(LI^2)$. In this case the succedent of the conclusion of $\mathcal{D}_2$ is $\Delta, B_{t_1}$, where $B_{t_1}$ is an atomic formula. There are two cases. 

\noindent (III.a) The principal formula $Ix[F, \exists ! x]$ of $(LI^2)$ is not one of the $A^k$: apply the induction hypothesis to the premises of $(LI^2)$ and then apply the rule.

\noindent (III.b) The principal formula $Ix[F, \exists ! x]$ of $(LI^2)$ is one of the $A^k$. Let $\Xi$ be $Ix[F, \exists !x]^{k-1}$, i.e. $\Xi$ consists of $k-1$ occurrences of $Ix[F, \exists !x]$, then $\mathcal{D}_2$ ends with: 

\bigskip

\begin{center}
\def\defaultHypSeparation{\hskip .075in}
\AxiomC{$\Xi, \Gamma\Rightarrow \Delta, F_{t_1}^x$}
\AxiomC{$\Xi, \Gamma\Rightarrow \Delta, F_{t_2}^x$}
\AxiomC{$\Xi, \Gamma\Rightarrow \Delta, \exists !t_1$}
\AxiomC{$\Xi, \Gamma\Rightarrow \Delta, \exists !t_2$}
\AxiomC{$\Xi, \Gamma\Rightarrow \Delta, B_{t_2}^x$}
\QuinaryInfC{$Ix[F, \exists !x]^k, \Gamma\Rightarrow \Delta, B_{t_1}^x$}
\DisplayProof
\end{center}

\bigskip

\noindent By induction hypothesis, we have:\bigskip

(1) $\Theta^{k-1}, \Gamma\Rightarrow \Lambda^{k-1}, \Delta, F_{t_1}^x$ 

(2) $\Theta^{k-1}, \Gamma\Rightarrow\Lambda^{k-1}, \Delta, F_{t_2}^x$

(3) $\Theta^{k-1}, \Gamma\Rightarrow \Lambda^{k-1},  \Delta, \exists !t_1$ 

(4) $\Theta^{k-1}, \Gamma\Rightarrow \Lambda^{k-1}, \Delta, \exists !t_2$

(5) $\Theta^{k-1}, \Gamma\Rightarrow \Lambda^{k-1}, \Delta, B_{t_2}^x$\bigskip

\noindent $A$ is principal in $\mathcal{D}_1$, so it ends with an application of $(RI)$ with $G$ replaced by $\exists !$:  

\bigskip

\begin{center} 
\AxiomC{$\Theta\Rightarrow \Lambda, F_t^x$}
\AxiomC{$\Theta\Rightarrow \Lambda, \exists ! t$}
\AxiomC{$\Theta\Rightarrow \Lambda, \exists !t$}
\AxiomC{$\exists ! a, F_a^x, \Theta \Rightarrow \Lambda, a=t$}
\QuaternaryInfC{$\Theta \Rightarrow \Lambda, Ix[F, \exists ! x]$}
\DisplayProof
\end{center} 

\bigskip

\noindent By the Substitution Lemma from the fourth premise:\bigskip

(6) $\exists ! t_1, F_{t_1}^x, \Theta \Rightarrow \Lambda, t_1=t$

(7) $\exists ! t_2, F_{t_2}^x, \Theta \Rightarrow \Lambda, t_2=t$\bigskip

\noindent From $B_{t_1}^x\Rightarrow  B_{t_1}^x$ by $(=E)$: (8) $t_1=t_2, B_{t_2}^x\Rightarrow  B_{t_1}^x$, and similarly (9) $t_1=t, t_2=t\Rightarrow t_1=t_2$. Two cuts and contraction with (6), (7) and (9) twice gives: (10) $\exists ! t_1, F_{t_1}^x, \exists ! t_2, F_{t_2}^x, \Theta \Rightarrow \Lambda, t_1=t_2$. A cut with (8) gives\bigskip

(11) $\exists ! t_1, F_{t_1}^x, \exists ! t_2, F_{t_2}^x, B_{t_2}^x, \Theta \Rightarrow \Lambda, B_{t_1}^x$\bigskip

\noindent Five cuts with (11) and (1)-(5) and contraction give $\Theta^k, \Gamma\Rightarrow \Lambda^k, \Delta, B_{t_1}$, which was to be proved.\bigskip

\noindent (IV) The last step of $\mathcal{D}_2$ is by $(LI^3)$. Two cases: 

\noindent (IV.a) The principal formula $Ix[F, \exists !x]$ of $(LI^3)$ is not one of the $A^k$: apply the induction hypothesis to the premises of $(LI^3)$ and then apply the rule. 

\noindent (IV.b) The principal formula $Ix[F, \exists !x]$ of $(LI^3)$ is one of the $A^k$. Then $\mathcal{D}_2$ ends with: 

\bigskip

\begin{center} 
\Axiom$F_a^x, \exists !_a^x, Ix[F, \exists !x]^{k-1}, \Gamma\fCenter \Delta$
\UnaryInf$Ix[F, \exists !x]^k, \Gamma\fCenter\Delta$
\DisplayProof
\end{center}

\bigskip

\noindent By induction hypothesis we have $ F_a^x, \exists !_a^x, \Theta^{k-1}, \Gamma\Rightarrow \Lambda^{k-1}, \Delta$, and so by the Substitution Lemma:\bigskip

(1) $F_t^x, \exists !_t^x, \Theta^{k-1}, \Gamma\Rightarrow \Lambda^{k-1}, \Delta$\bigskip

\noindent $A$ is principal in $\mathcal{D}_1$, so it ends with an application of $(RI)$ with $G$ replaced by $\exists !$: 

\bigskip

\begin{center} 
\AxiomC{$\Theta\Rightarrow \Lambda, F_t^x$}
\AxiomC{$\Theta\Rightarrow \Lambda, \exists ! t$}
\AxiomC{$\Theta\Rightarrow \Lambda, \exists !t$}
\AxiomC{$\exists ! a, F_a^x, \Theta \Rightarrow \Lambda, a=t$}
\QuaternaryInfC{$\Theta \Rightarrow \Lambda, Ix[F, \exists ! x]$}
\DisplayProof
\end{center} 

\bigskip

\noindent Apply two cuts with the first, and second or third, premise and (1), and contraction to conclude $\Theta^k, \Gamma\Rightarrow \Lambda^k, \Delta$.\bigskip

\noindent (V) The last step of $\mathcal{D}_2$ is by $(LI^4)$. As in case (III), the succedent of the conclusion of $\mathcal{D}_2$ is $\Delta, B_{t_1}$, where $B_{t_1}$ is an atomic formula. Two cases: 

\noindent (V.a) The principal formula $Ix[F, x=t_2]$ of $(LI^4)$ is not one of the $A^k$: apply the induction hypothesis to the premises of $(LI^3)$ and then apply the rule. 

\noindent (V.b) The principal formula $Ix[F, x=t_2]$ of $(LI^4)$ is one of the $A^k$. Let $\Xi$ be $Ix[F, x=t_2]^{k-1}$, i.e. $\Xi$ consists of $k-1$ occurrences of $Ix[F, x=t_2]$, then $\mathcal{D}_2$ ends with:

\bigskip

\begin{center}
\AxiomC{$\Xi, \Gamma\Rightarrow \Delta, F_{t_1}^x$}
\AxiomC{$\Xi, \Gamma\Rightarrow \Delta, \exists !t_1$}
\AxiomC{$\Xi, \Gamma\Rightarrow \Delta, \exists !t_2$}
\AxiomC{$\Xi, \Gamma\Rightarrow \Delta, B_{t_2}^x$}
\QuaternaryInfC{$Ix[F, x=t_2]^k, \Gamma\Rightarrow \Delta, B_{t_1}^x$}
\DisplayProof
\end{center}

\bigskip

\noindent By induction hypothesis, we have the following, although we won't need (3):\bigskip

(1) $\Theta^{k-1}, \Gamma\Rightarrow \Lambda^{k-1}, \Delta, F_{t_1}^x$

(2) $\Theta^{k-1}, \Gamma\Rightarrow \Lambda^{k-1}, \Delta, \exists !t_1$

(3) $\Theta^{k-1}, \Gamma\Rightarrow \Lambda^{k-1}, \Delta, \exists !t_2$

(4) $\Theta^{k-1}, \Gamma\Rightarrow \Lambda^{k-1}, \Delta, B_{t_2}^x$\bigskip

\noindent $A$ is principal in $\mathcal{D}_1$, so it ends with an application of $(RI)$ with $Gx$ replaced by $x=t_2$:

\bigskip

\begin{center} 
\AxiomC{$\Theta\Rightarrow \Lambda, F_t^x$}
\AxiomC{$\Theta\Rightarrow \Lambda, t=t_2$}
\AxiomC{$\Theta\Rightarrow \Lambda, \exists !t$}
\AxiomC{$\exists ! a, F_a^x, \Theta \Rightarrow \Lambda, a=t$}
\QuaternaryInfC{$\Theta \Rightarrow \Lambda, Ix[F, x=t_2]$}
\DisplayProof
\end{center}

\bigskip

\noindent As in case (III.b), we have (5) $t_1=t_2, B_{t_2}^x\Rightarrow  B_{t_1}^x$ and (6) $t_1=t, t_2=t\Rightarrow t_1=t_2$, and from (6) and the second premise by cut: $t_1=t, \Theta\Rightarrow \Lambda, t_1=t_2$, from which by (5) and (4) by two times cut:\bigskip

(7) $t_1=t, \Theta^k, \Gamma\Rightarrow \Lambda^k, \Delta, B_{t_1}^x$\bigskip 

\noindent By the Substitution Lemma from the fourth premise of the final $(RI)$ of $\mathcal{D}_1$:\bigskip

(8) $\exists ! t_1, F_{t_1}^x, \Theta \Rightarrow \Lambda, t_1=t$\bigskip

\noindent whence from (7) by cut and contraction, $\exists ! t_1, F_{t_1}^x, \Theta^k, \Gamma \Rightarrow \Lambda^k, \Delta, B_{t_1}^x$, from which by cut and contraction with (1) and (2) (or also the first and second premise of $(RI)$) $\Theta^k, \Gamma \Rightarrow \Lambda^k, \Delta, B_{t_1}^x$, which was to be proved.\bigskip

\noindent (VI) The last step of $\mathcal{D}_2$ is by $(LI^5)$. Two cases: 

\noindent (VI.a) The principal formula $Ix[F, x=t]$ of $(LI^5)$ is not one of the $A^k$: apply the induction hypothesis to the premises of $(LI^5)$ and then apply the rule. 

\noindent (VI.b) The principal formula $Ix[F, x=t]$ of $(LI^5)$ is one of the $A^k$. Then $\mathcal{D}_2$ ends with:

\bigskip

\begin{center}
\AxiomC{$Ix[F, x=t]^{k-1}, \Gamma\Rightarrow \Delta, \exists !t$}
\AxiomC{$F_a^x, \exists !a, Ix[F, x=t]^{k-1}, \Gamma\Rightarrow \Delta$}
\BinaryInfC{$Ix[F, x=t]^k, \Gamma\Rightarrow\Delta$}
\DisplayProof
\end{center} 

\bigskip

\noindent By induction hypothesis, we have the following, of which we need only the second:\bigskip

(1) $\Theta^{k-1}, \Gamma\Rightarrow \Lambda^{k-1}, \Delta, \exists !t$

(2) $F_a^x, \exists !a, \Theta^{k-1}, \Gamma\Rightarrow \Lambda^{k-1}, \Delta$\bigskip

\noindent $A$ is principal in $\mathcal{D}_1$, so it ends with an application of $(RI)$ with $Gx$ replaced by $x=t$:

\bigskip

\begin{center} 
\AxiomC{$\Theta\Rightarrow \Lambda, F_{t_1}^x$}
\AxiomC{$\Theta\Rightarrow \Lambda, t_1=t$}
\AxiomC{$\Theta\Rightarrow \Lambda, \exists !t_1$}
\AxiomC{$\exists ! a, F_a^x, \Theta \Rightarrow \Lambda, a=t_1$}
\QuaternaryInfC{$\Theta \Rightarrow \Lambda, Ix[F, x=t]$}
\DisplayProof
\end{center}

\bigskip

\noindent By the Substitution Lemma from (2):\bigskip 

(3) $F_{t_1}^x, \exists !t_1, \Theta^{k-1}, \Gamma\Rightarrow \Lambda^{k-1}, \Delta$\bigskip

\noindent from which by the first and third premises of $(RI)$ with cut and contraction $\Theta^k, \Gamma\Rightarrow \Lambda^k, \Delta$.\bigskip

\noindent This completes the proof of the Right Reduction Lemma.\bigskip 

\noindent A note on steps (III.b) and (V.b) might be in order: cuts on identities are eliminable from GPFL$_=$. This does not change in GPFL$_=^I$, as identity is not principal in any of its rules. Notice incidentally that this would not be the case had we chosen $(LI{^2}')$ (or a corresponding version of $(LI^4)$) instead of $(LI^{2})$ (and $(LI^4)$. Thus the cuts on identities appealed to in steps (III.b) and (V.b) are eliminable and they are of course, being cuts on atomic formulas, of lower degree than $Ix[F, G]$, $Ix[F, \exists!x]$ and $Ix[F, x=t]$. Recall again also that the formula $B$ in $(LI^2)$ and $(LI^4)$ is atomic. Thus all cuts applied in the proof are of lower degree than the formula $A$ of the conclusions of $\mathcal{D}_1$ and $\mathcal{D}_2$. 

\bigskip

\begin{lemma}[Left Reduction]
If $\mathcal{D}_1\vdash \Gamma\Rightarrow \Delta, A^k$ and $\mathcal{D}_2\vdash A, \Theta\Rightarrow \Lambda$ have degrees $d(\mathcal{D}_1), d(\mathcal{D}_2)<d(A)$, then there is a proof $\mathcal{D}\vdash \Gamma, \Theta^k\Rightarrow \Delta, \Lambda^k$ with $d(\mathcal{D})<d(A)$. 
\end{lemma} 

\noindent \emph{Proof} by induction over the height of $\mathcal{D}_1$. 

The basis is trivial, as then $\mathcal{D}_1$ is an axiom, and $\Gamma$ consists of one occurrence of $A$ and $\Delta$ is empty. What needs to be shown is that $A, \Theta\Rightarrow \Lambda$, which is already given by $\mathcal{D}_2$.

For the induction step, we distinguish two cases, and again we continue Indrzejczak's proof by adding the new cases arising in GPFL$_=^I$ through the addition of $I$.\bigskip

\noindent (A) None of the $A^k$ in the succedent of the conclusion of $\mathcal{D}_1$ is principal. Then we apply the induction hypothesis to the premises of the final rule applied in $\mathcal{D}_1$ and apply the final rule once more.\bigskip 

\noindent (B) Some of the $A^k$ in the succedent of the conclusion of $\mathcal{D}_1$ are principal. Then there are three options.\bigskip 

\noindent (I) The final rule applied in $\mathcal{D}_1$ is $(RI)$. Let $\Xi$ be $Ix[F, G]^{k-1}$, i.e. $\Xi$ consists of $k-1$ occurrences of $Ix[F, G]$, then $\mathcal{D}_1$ ends with:

\bigskip 

\begin{center}
\AxiomC{$\Gamma\Rightarrow \Delta, \Xi, F_t^x$}
\AxiomC{$\Gamma\Rightarrow \Delta, \Xi, G_t^x$}
\AxiomC{$\Gamma\Rightarrow \Delta, \Xi, \exists !t$}
\AxiomC{$\exists ! a, F_a^x, \Gamma \Rightarrow \Delta, \Xi, a=t$}
\QuaternaryInfC{$\Gamma\Rightarrow \Delta, Ix[F, G]^k$}
\DisplayProof
\end{center}

\bigskip

\noindent By induction hypothesis, we have\bigskip

(1) $\Gamma, \Theta^{k-1} \Rightarrow \Delta, \Lambda^{k-1}, F_t^x$

(2) $\Gamma, \Theta^{k-1}\Rightarrow \Delta, \Lambda^{k-1}, G_t^x$

(3)  $\Gamma, \Theta^{k-1}\Rightarrow \Delta, \Lambda^{k-1}, \exists !t$

(4) $\exists !a, F_a^x, \Gamma, \Theta^{k-1}\Rightarrow \Delta, \Lambda^{k-1}, a=t$\bigskip

\noindent Apply $(RI)$ with (1) to (4) as premises to conclude\bigskip 

(5) $\Gamma, \Theta^{k-1} \Rightarrow \Delta, \Lambda^{k-1}, Ix[F, G]$\bigskip 

\noindent Here $Ix[F, G]$ is principal, so we apply the Right Reduction Lemma to the deduction concluding (5) and $\mathcal{D}_2$ (where $k=1$) to conclude $\Gamma, \Theta^k \Rightarrow \Delta, \Lambda^k$. \bigskip

\noindent (II) The final rule applied in $\mathcal{D}_1$ is $(LI^2)$. Let $\Xi$ be ${B_{t_1}^x}^{k-1}$, i.e. $\Xi$ consists of $k-1$ occurrences of ${B_{t_1}^x}^{k-1}$, then $\mathcal{D}_1$ ends with:

\bigskip 

\begin{center}
\def\defaultHypSeparation{\hskip .05in}
\AxiomC{$\Gamma\Rightarrow \Delta, \Xi, F_{t_1}^x$}
\AxiomC{$\Gamma\Rightarrow \Delta, \Xi, F_{t_2}^x$}
\AxiomC{$\Gamma\Rightarrow \Delta, \Xi, \exists !t_1$}
\AxiomC{$\Gamma\Rightarrow \Delta, \Xi, \exists !t_2$}
\AxiomC{$\Gamma\Rightarrow \Delta, \Xi, B_{t_2}^x$}
\QuinaryInfC{$Ix[F, \exists !x],\Gamma\Rightarrow \Delta,  {B_{t_1}^x}^k$}
\DisplayProof
\end{center}

\bigskip

\noindent By induction hypothesis, we have: \bigskip

(1) $\Gamma, \Theta^{k-1} \Rightarrow \Delta, \Lambda^{k-1}, F_{t_1}^x$

(2) $\Gamma, \Theta^{k-1} \Rightarrow \Delta, \Lambda^{k-1}, F_{t_2}^x$

(3) $\Gamma, \Theta^{k-1} \Rightarrow \Delta, \Lambda^{k-1}, \exists !t_1$

(4) $\Gamma, \Theta^{k-1} \Rightarrow \Delta, \Lambda^{k-1}, \exists !t_2$

(5) $\Gamma, \Theta^{k-1} \Rightarrow \Delta, \Lambda^{k-1}, B_{t_2}^x$\bigskip

\noindent Apply $(LI^2)$ with (1) to (5) as premises to conclude\bigskip

(6) $Ix[F, \exists !x], \Gamma, \Theta^{k-1} \Rightarrow \Delta, \Lambda^{k-1}, B_{t_1}^x$\bigskip 

\noindent Again $B_{t_1}^x$ is principal, so we apply the Right Reduction Lemma to the deduction concluding (6) and $\mathcal{D}_2$ (where $k=1$) to conclude $Ix[F, \exists !x], \Gamma, \Theta^k \Rightarrow \Delta, \Lambda^k$. \bigskip

\noindent (III) The final rule applied in $\mathcal{D}_1$ is $(LI^4)$. Let $\Xi$ be ${B_{t_1}^x}^{k-1}$, i.e. $\Xi$ consists of $k-1$ occurrences of ${B_{t_1}^x}^{k-1}$, then $\mathcal{D}_1$ ends with:

\bigskip

\begin{center}
\AxiomC{$\Gamma\Rightarrow \Delta, \Xi, F_{t_1}^x$}
\AxiomC{$\Gamma\Rightarrow \Delta, \Xi, \exists !t_1$}
\AxiomC{$\Gamma\Rightarrow \Delta, \Xi, \exists !t_2$}
\AxiomC{$\Gamma\Rightarrow \Delta, \Xi, B_{t_2}^x$}
\QuaternaryInfC{$Ix[F, x=t_2], \Gamma\Rightarrow \Delta, {B_{t_1}^x}^k$}
\DisplayProof
\end{center}

\bigskip

\noindent By induction hypothesis, we have\bigskip

(1) $\Gamma, \Theta^{k-1} \Rightarrow \Delta, \Lambda^{k-1}, F_{t_1}^x$

(2) $\Gamma, \Theta^{k-1}\Rightarrow \Delta, \Lambda^{k-1}, \exists !t_1$

(3) $\Gamma, \Theta^{k-1} \Rightarrow \Delta, \Lambda^{k-1}, \exists !t_2$

(4) $\Gamma, \Theta^{k-1} \Rightarrow \Delta, \Lambda^{k-1}, B_{t_2}^x$\bigskip

\noindent Apply $(LI^4)$ with (1) to (4) as premises to conclude\bigskip

(5) $Ix[F, x=t_2], \Theta^{k-1} \Gamma\Rightarrow \Delta, \Lambda^{k-1}, {B_{t_1}^x}^k$\bigskip

\noindent Once more $B_{t_1}^x$ is principal, so apply the Right Reduction Lemma to the deduction concluding (5) and $\mathcal{D}_2$ (where $k=1$) to conclude $Ix[F, x=t_2], \Gamma, \Theta^k \Rightarrow \Delta, \Lambda^k$. \bigskip

\noindent This completes the proof of the Left Reduction Lemma. 

\bigskip

\begin{theorem}[Cut Elimination] 
For every deduction in GPFL$_=^I$, there is a deduction that is free of cuts. 
\end{theorem}

\noindent \emph{Proof}. The theorem follows from the Right and Left Reduction Lemmas by induction over the degree of the proof, with subsidiary deductions over the number of cut formulas of highest degree, as in Indrzejczak's paper.

\section{Concluding Remarks} 
\subsection{Slightly Simpler Rules for $I$}
It is possible to simplify the rules for $I$ in the sense that most rules can be put into a form that requires fewer premises by putting existence assumptions into the antecedents of the conclusions instead of the consequents premises:\footnote{This possibility was pointed out by a referee for \emph{Studia Logica}, to whom many thanks for the suggestion.} 

\bigskip

\noindent
\AxiomC{$\Gamma\Rightarrow \Delta, F_t^x$}
\AxiomC{$\Gamma\Rightarrow \Delta, G_t^x$}
\AxiomC{$\exists ! a, F_a^x, \Gamma \Rightarrow \Delta, a=t$}
\LeftLabel{$(RI^S)$ \ }
\TrinaryInfC{$\exists ! t, \Gamma \Rightarrow \Delta, I x[F, G]$}
\DisplayProof

\bigskip

\noindent where $a$ does not occur in the conclusion. 

\bigskip

\noindent
\AxiomC{$\Gamma\Rightarrow \Delta, F_t^x$}
\AxiomC{$F_a^x, \exists ! a, \Gamma\Rightarrow \Delta, a=t$}
\AxiomC{$F_b^x, G_b^x, \exists !b, \Gamma\Rightarrow \Delta$}
\LeftLabel{$(LI^{1S})$ \ }
\TrinaryInfC{$\exists !t, Ix[F, G], \Gamma\Rightarrow \Delta$}
\DisplayProof

\bigskip

\noindent where $a$ and $b$ do not occur in the conclusion. 

\bigskip

\noindent
\AxiomC{$\Gamma\Rightarrow \Delta, F_{t_1}^x$}
\AxiomC{$\Gamma\Rightarrow \Delta, F_{t_2}^x$}
\AxiomC{$\Gamma\Rightarrow \Delta, A_{t_2}^x$}
\LeftLabel{$(LI^{2S})$ \ }
\TrinaryInfC{$\exists ! t_1, \exists !t_2, Ix[F, \exists !x], \Gamma\Rightarrow \Delta, A_{t_1}^x$}
\DisplayProof

\bigskip

\noindent where $A$ is an atomic formula.

\bigskip

\noindent
\Axiom$F_a^x, \exists !a, \Gamma\fCenter \Delta$
\LeftLabel{$(LI^3)$ \ }
\UnaryInf$Ix[F, \exists !x], \Gamma\fCenter\Delta$
\DisplayProof

\bigskip

\noindent where $a$ does not occur in the conclusion.

\bigskip

\noindent
\AxiomC{$\Gamma\Rightarrow \Delta, F_{t_1}^x$}
\AxiomC{$\Gamma\Rightarrow \Delta, A_{t_2}^x$}
\LeftLabel{$(LI^{4S})$ \ }
\BinaryInfC{$\exists !t_1, \exists t_2, Ix[F, x=t_2], \Gamma\Rightarrow \Delta, A_{t_1}^x$}
\DisplayProof

\bigskip

\noindent where $A$ is an atomic formula.

\bigskip

\noindent
\AxiomC{$F_a^x, \exists !a, \Gamma\Rightarrow \Delta$}
\LeftLabel{$(LI^{5S})$ \ }
\UnaryInfC{$\exists !t, Ix[F, x=t], \Gamma\Rightarrow\Delta$}
\DisplayProof

\bigskip

\noindent where $a$ does not occur in the conclusion. 

\bigskip

\noindent This reduces the branching factor of deductions, which helps with proof search. The rules of Section 3 have the advantage of corresponding a little more directly to the rules of natural deduction give in a previous paper, and in every rule there is at most one principal formula to the left or to the right of $\Rightarrow$ in the conclusion. 

If we go for the simplified rules of this section, it makes sense to change $(LI^{2S})$ and $(LI^{4S})$ in the way proposed by Indrzejczak (cf. footnote \ref{andrzejfootnote}) to: 

\bigskip

\noindent 
\AxiomC{$\Gamma\Rightarrow \Delta, F_{t_1}^x$}
\AxiomC{$\Gamma\Rightarrow \Delta, F_{t_2}^x$}
\AxiomC{$t_!=t_2, \Gamma\Rightarrow \Delta$}
\LeftLabel{$(LI^{2SI})$ \ }
\TrinaryInfC{$\exists ! t_1, \exists !t_2, Ix[F, \exists !x], \Gamma\Rightarrow \Delta$}
\DisplayProof

\bigskip

\noindent
\AxiomC{$\Gamma\Rightarrow \Delta, F_{t_1}^x$}
\AxiomC{$t_1=t_2, \Gamma\Rightarrow \Delta$}
\LeftLabel{$(LI^{4SI})$ \ }
\BinaryInfC{$\exists !t_1, \exists t_2, Ix[F, x=t_2], \Gamma\Rightarrow \Delta$}
\DisplayProof

\bigskip

\noindent In the system consisting of $(RI)$, $(LI^{1S})$,  $(LI^{2SI})$,  $(LI^{3S})$,  $(LI^{4SI})$ and $(LI^{5S})$ added to GPFL$_=$, of the  rules for $I$ only $(RI)$ introduces a principal formula to the right of $\Rightarrow$, and formulas of the form $\exists!t $ are never principal in that position. Thus steps (I), (II.a), (III.a), (IV.a), (V.a) and (VI.a) of the Right Reduction Lemma go through as before, with some minor rephrasing (as $A$ is principal in the final step of $\mathcal{D}_1$, the new existence formulas to the left of $\Rightarrow$ cannot be amongst the $A^k$). 

Step (II.b) goes through almost as before. We now have an additional $\exists !t$ in the antecedent of the conclusion of $\mathcal{D}_2$. We still have sequent (5), which contains the required $\exists !t$, so we apply Cut to the first two premises of the application of $(RI^S)$ with which $\mathcal{D}_1$ ends, and we're done. 

Step (III.b) also goes through almost as before and along a similar pattern as new case (II.b). We now have additional $\exists !t_1$ and $\exists !t_2$ in the antecedent of the conclusion of $\mathcal{D}_2$. We still have sequents (6) and (7), which contain the required $\exists !t_1$ and $\exists !t_2$, and we apply Cut twice to them and $t_1=t, t_2=t\Rightarrow t_1=t_2$, then to the resulting sequent and the third premise of the application of $(LI^{2SI})$ with which $\mathcal{D}_1$ ends, and we're done. 

Step (IV.b) concerns rule $(LI^3)$ which remains unchanged. Step (V.b) goes through with adjustments similar to those in the adjusted step (III.b); step (VI.b) similar to adjusted step (II.b). 

In the Left Reduction Lemma, we only need to consider the case where the final rule applied in $\mathcal{D}_1$ is $(RI^S)$, and as in previous cases, we still have sequent (5) so all is well.

\subsection{$I$ in Negative Free and Classical Logic}
Indrzejczak's system GNFL$_=$ of negative free logic arises from GPFL$_=$ by changing $(=E)$ to the rule $(NEI)$ below and adding the \emph{rules of strictness}:

\bigskip

\noindent
\Axiom$t=t, \Gamma \fCenter \Delta$
\LeftLabel{$(NEI) \ $}
\UnaryInf$\exists ! t, \Gamma \fCenter\Delta$
\DisplayProof

\bigskip

\noindent
\Axiom$\exists ! t_i, \Gamma\fCenter \Delta$
\LeftLabel{$(NEE) \ $}
\UnaryInf$Rt_1\ldots t_n, \Gamma \fCenter \Delta$
\DisplayProof\qquad
\Axiom$\exists ! t_i, \Gamma\fCenter \Delta$
\LeftLabel{$(NEE') \ $}
\UnaryInf$\exists ! ft_1\ldots t_n, \Gamma \fCenter \Delta$
\DisplayProof

\bigskip

\noindent for $i\leq n$, for all $n$-place predicates $R$ and functions $f$. 

The following are appropriate rules for the binary quantifier $I$ in negative free logic:

\bigskip

\noindent
\AxiomC{$\Gamma\Rightarrow \Delta, F_t^x$}
\AxiomC{$\Gamma\Rightarrow \Delta, G_t^x$}
\AxiomC{$\Gamma\Rightarrow \Delta, \exists !t$}
\AxiomC{$F_a^x, \Gamma \Rightarrow \Delta, a=t$}
\LeftLabel{$(RI^N)$ \ }
\QuaternaryInfC{$\Gamma \Rightarrow \Delta, Ix[F, G]$}
\DisplayProof

\bigskip

\noindent \
\Axiom$F_a^x, G_a^x, \exists !a, \Gamma \fCenter \Delta$
\LeftLabel{$(LI^{N1})$ \ }
\UnaryInf$Ix[F, G], \Gamma\fCenter\Delta$
\DisplayProof

\bigskip

\noindent
\def\defaultHypSeparation{\hskip .1in}
\AxiomC{$\Gamma \Rightarrow \Delta, \exists !t_1$}
\AxiomC{$\Gamma \Rightarrow \Delta, \exists !t_2$}
\AxiomC{$\Gamma \Rightarrow \Delta, F_{t_1}^x$}
\AxiomC{$\Gamma \Rightarrow \Delta, F_{t_2}^x$}
\AxiomC{$\Gamma \Rightarrow \Delta, A_{t_2}^x$}
\LeftLabel{$(LI^{N2})$ \ }
\QuinaryInfC{$Ix[F, G], \Gamma\Rightarrow \Delta, A_{t_1}$}
\DisplayProof

\bigskip

\noindent where in $(RI^N)$ and $(LI^{N1})$, $a$ does not occur in the conclusion, and in $(LI^{N2})$ $A$ is an atomic formula.\footnote{These are the rules of \cite{kurbisiotaI} transposed to sequent calculus. They could be simplified, analogously to the proposal of the previous section, by deleting the second premise of $(RI^N)$ and adding $\exists !t$ to the antecedent of the conclusion, and by deleting the first two premises of $(LI^{N2})$ and adding $\exists !t_1$ and $\exists !t_2$ to the antecedent of the conclusion.}\bigskip

\noindent Let GNFL$_=^I$ be GNFL$_=$ with its language extended by $I$ and $(RI^N)$, $(LI^{N1})$ and $(LI^{N2})$ added as rules. In this system $\vdash Ix[F, G]\leftrightarrow \exists x(\forall y(F_y^x\leftrightarrow x=y)\land G_y^x)$. Thus it is adequate as a formalisation of a Russellian theory of definite descriptions with scope distinctions marked by the square brackets of the binary quantifier $I$. Cut elimination is provable for GNFL$_=^I$, too, following once more Indrzejczak's proof of cut elimination for GNFL$_=$ and extending it by the new cases for $I$. 

Finally, one could even consider adding $I$ to classical logic. For that we would need to change the rules for the quantifiers of GPFL$_=$ in well known fashion, and then suitable rules for $I$ are the following:\bigskip

\noindent
\AxiomC{$\Gamma\Rightarrow \Delta, F_t^x$}
\AxiomC{$\Gamma\Rightarrow \Delta, G_t^x$}
\AxiomC{$F_a^x, \Gamma \Rightarrow \Delta, a=t$}
\LeftLabel{$(RI^C)$ \ }
\TrinaryInfC{$\Gamma \Rightarrow \Delta, Ix[F, G]$}
\DisplayProof

\bigskip

\noindent \
\Axiom$F_a^x, G_a^x, \Gamma \fCenter \Delta$
\LeftLabel{$(LI^{C1})$ \ }
\UnaryInf$Ix[F, G], \Gamma\fCenter\Delta$
\DisplayProof

\bigskip

\noindent
\AxiomC{$\Gamma\Rightarrow \Delta, F_{t_1}^x$}
\AxiomC{$\Gamma\Rightarrow \Delta, F_{t_2}^x$}
\AxiomC{$\Gamma\Rightarrow \Delta, A_{t_2}^x$}
\LeftLabel{$(LI^{C2})$ \ }
\TrinaryInfC{$Ix[F, G], \Gamma \Rightarrow \Delta, A_{t_1}$}
\DisplayProof

\bigskip

\noindent where in $(RI^C)$ and $(LI^{C1})$, $a$ does not occur in the conclusion, and in $(LI^{C2})$ $A$ is an atomic formula.\bigskip

\noindent Then $Ix[F, G]\leftrightarrow \exists x(\forall y(F_y^x\leftrightarrow x=y)\land G_y^x)$ is also provable, and Cut elimination goes through as before. 

It is interesting to note that the rules for $I$ in negative free logic and classical logic are significantly simpler than those for $I$ in positive free logic. The reasons is that in both the former logics, $Ix[F, G]$ is equivalent to a formula that, albeit already fairly complex, is still reasonably straightforward, namely the formula that expresses the Russellian analysis of `The $F$ is $G$'. Thus $Ix[F, G]$ is definable in terms or eliminable in favour of the latter and all we are required to do, should we wish to keep it as a primitive nonetheless, is to pretend to look for rules that would allow to introduce $\exists x(\forall y(F_y^x\leftrightarrow x=y)\land G_y^x)$ immediately to the left and to right of the sequent arrow, and then use those rules for $I$ instead. The situation is more complicated in positive free logic, as there $Ix[F, G]$ is not straightforwardly equivalent to anything else: adding a means for formalising definite descriptions to positive free logic constitutes a genuine extension of its expressive power. $Ix[F, G]$ is equivalent to $\exists x(\forall y(F_y^x\leftrightarrow x=y)\land G_y^x)$ only under the assumption that a unique $F$ exists. The latter is also already fairly complex; indeed, it is expressible by a formula involving $I$. `The $F$ is $G$' says something rather more intricate in positive free logic than it does in negative free logic. This is again to do with the aim of theorist of definite descriptions who prefer positive free logic, discussed in the introduction, to avoid commitment to the existence of a unique $F$ with an assertion of `The $F$ is $G$', and, indeed, to commit to nothing much at all should there not be one.

\bigskip

\setlength{\bibsep}{0pt}
\bibliographystyle{chicago}
\bibliography{KurbisIotaSequentI}

\begin{thebibliography}{}

\bibitem[\protect\citeauthoryear{Bencivenga}{Bencivenga}{1986}]{bencivengahandbook}
Bencivenga, E. (1986).
\newblock Free logics.
\newblock In D.~Gabbay and F.~Guenther (Eds.), {\em Handbook of Philosophical
  Logic. Volume III: Alternatives to Classical Logic}, pp.\  373--426.
  Dortrecht: Springer.

\bibitem[\protect\citeauthoryear{Bostock}{Bostock}{1997}]{bostockintermediate}
Bostock, D. (1997).
\newblock {\em Intermediate Logic}.
\newblock Oxford: Clarendon Press.

\bibitem[\protect\citeauthoryear{Czermak}{Czermak}{1974}]{czermakdefdescr}
Czermak, J. (1974).
\newblock A logical calculus with definite descriptions.
\newblock {\em Journal of Philosophical Logic\/}~{\em 3\/}(3), 211--228.

\bibitem[\protect\citeauthoryear{Dummett}{Dummett}{1981}]{dummettfregelanguage}
Dummett, M. (1981).
\newblock {\em Frege. Philosophy of Language\/} (2 ed.).
\newblock London: Duckworth.

\bibitem[\protect\citeauthoryear{Fitting and Mendelsohn}{Fitting and
  Mendelsohn}{1998}]{mendelsohnfitting}
Fitting, M. and R.~L. Mendelsohn (1998).
\newblock {\em First-Order Modal Logic}.
\newblock Dordrecht, Boston, London: Kluwer.

\bibitem[\protect\citeauthoryear{Garson}{Garson}{2013}]{garsonmodallogic}
Garson, J.~W. (2013).
\newblock {\em Modal Logic for Philosophers\/} (2 ed.).
\newblock Cambridge University Press.

\bibitem[\protect\citeauthoryear{Gratzl}{Gratzl}{2015}]{gratzldefdescr}
Gratzl, N. (2015).
\newblock Incomplete symbols -- definite descriptions revisited.
\newblock {\em Journal of Philosophical Logic\/}~{\em 44\/}(5), 489--506.

\bibitem[\protect\citeauthoryear{Hintikka}{Hintikka}{1959}]{hintikkatowardsdefdesc}
Hintikka, J. (1959).
\newblock Towards a theory of definite descriptions.
\newblock {\em Analysis\/}~{\em 19\/}(4), 79--85.

\bibitem[\protect\citeauthoryear{Indrzejczak}{Indrzejczak}{2018a}]{andrzejmodaldescription}
Indrzejczak, A. (2018a).
\newblock Cut-free modal theory of definite descriptions.
\newblock In G.~M. G.~Bezhanishvili, G.~D'Agostino and T.~Studer (Eds.), {\em
  Advances in Modal Logic}, Volume~12, pp.\  359--378. London: College
  Publications.

\bibitem[\protect\citeauthoryear{Indrzejczak}{Indrzejczak}{2018b}]{andrzejfregean}
Indrzejczak, A. (2018b).
\newblock Fregean description theory in proof-theoretical setting.
\newblock {\em Logic and Logical Philosophy\/}~{\em 28\/}(1), 137--155.

\bibitem[\protect\citeauthoryear{Indrzejczak}{Indrzejczak}{2020a}]{andrzejexistencedefinedness}
Indrzejczak, A. (2020a).
\newblock Existence, definedness and definite descriptions in hybrid modal
  logic.
\newblock In N.~Olivetti, R.~Verbrugge, S.~Negri, and G.~Sandu (Eds.), {\em
  Advances in Modal Logic 13}. Rickmansworth: College Publications.

\bibitem[\protect\citeauthoryear{Indrzejczak}{Indrzejczak}{2020b}]{andrzejfreedefdescr}
Indrzejczak, A. (2020b).
\newblock Free definite description theory - sequent calculi and cut
  elimination.
\newblock {\em Logic and Logical Philosophy\/}~{\em 29\/}(4), Volume 29 (2020),
  505--539.

\bibitem[\protect\citeauthoryear{Indrzejczak}{Indrzejczak}{2021}]{andrzejcutfreefreelogic}
Indrzejczak, A. (2021).
\newblock Free logics are cut free.
\newblock {\em Studia Logica\/}~{\em online first}.

\bibitem[\protect\citeauthoryear{K\"urbis}{K\"urbis}{2019a}]{kurbisiotaI}
K\"urbis, N. (2019a).
\newblock A binary quantifier for definite descriptions in intuitionist
  negative free logic: Natural deduction and normalisation.
\newblock {\em Bulletin of the Section of Logic\/}~{\em 48\/}(2), 81--97.

\bibitem[\protect\citeauthoryear{K\"urbis}{K\"urbis}{2019b}]{kurbisiotaII}
K\"urbis, N. (2019b).
\newblock Two treatments of definite descriptions in intuitionist negative free
  logic.
\newblock {\em Bulletin of the Section of Logic\/}~{\em 48\/}(4), 299--318.

\bibitem[\protect\citeauthoryear{K\"urbis}{K\"urbis}{2021}]{kurbisiotaIII}
K\"urbis, N. (2021).
\newblock Definite descriptions in intuitionist positive free logic.
\newblock {\em Logic and Logical Philosophy\/}~{\em 30\/}(2), 327--358.

\bibitem[\protect\citeauthoryear{Lambert}{Lambert}{1961}]{lambertnotesEII}
Lambert, K. (1961).
\newblock {Notes on ''E!'': II}.
\newblock {\em Philosophical Studies\/}~{\em 12\/}(1/2), 1--5.

\bibitem[\protect\citeauthoryear{Lambert}{Lambert}{1962}]{lambertnotesEIII}
Lambert, K. (1962).
\newblock {Notes on ''E!'' III: A} theory of descriptions.
\newblock {\em Philosophical Studies\/}~{\em 13\/}(4), 51--59.

\bibitem[\protect\citeauthoryear{Lambert}{Lambert}{1964}]{lambertnotesEIV}
Lambert, K. (1964).
\newblock {Notes on ''E!'' IV: A} reduction in free quantification theory with
  identity and descriptions.
\newblock {\em Philosophical Studies\/}~{\em 15\/}(5), 85--88.

\bibitem[\protect\citeauthoryear{Lambert}{Lambert}{2001}]{lambertfreedef}
Lambert, K. (2001).
\newblock Free logic and definite descriptions.
\newblock In E.~Morscher and A.~Hieke (Eds.), {\em New Essays in Free Logic in
  Honour of Karel Lambert}. Dordrecht: Kluwer.

\bibitem[\protect\citeauthoryear{Morscher and Simons}{Morscher and
  Simons}{2001}]{morschersimonsfreelogic}
Morscher, E. and P.~Simons (2001).
\newblock Free logic: A fifty-year past and an open future.
\newblock In E.~Morscher and A.~Hieke (Eds.), {\em New Essays in Free Logic in
  Honour of Karel Lambert}. Dortrecht: Kluwer.

\bibitem[\protect\citeauthoryear{Neale}{Neale}{1990}]{nealedescriptions}
Neale, S. (1990).
\newblock {\em Descriptions}.
\newblock Cambridge, Mass.: MIT Press.

\bibitem[\protect\citeauthoryear{Ramsey}{Ramsey}{1990}]{ramseyphilosophy}
Ramsey, F.~P. (1990).
\newblock Philosophy.
\newblock In H.~Mellor (Ed.), {\em Philosophical Papers}. Cambridge University
  Press.

\bibitem[\protect\citeauthoryear{Russell and Whitehead}{Russell and
  Whitehead}{1910}]{russellwhitehead}
Russell, B. and A.~N. Whitehead (1910).
\newblock {\em Principia Mathematica}, Volume~1.
\newblock Cambridge University Press.

\bibitem[\protect\citeauthoryear{Troestra and Schwichtenberg}{Troestra and
  Schwichtenberg}{2000}]{troelstraschwichtenberg}
Troestra, A. and H.~Schwichtenberg (2000).
\newblock {\em Basic Proof Theory\/} (2 ed.).
\newblock Cambridge University Press.

\bibitem[\protect\citeauthoryear{van Fraassen}{van
  Fraassen}{1991}]{fraassenthexxlambert}
van Fraassen, B.~C. (1991).
\newblock On (the x) (x = {Lambert}).
\newblock In B.~S. Wolfgang~Spohn, Bas C. van~Fraassen (Ed.), {\em Existence
  and Explanation. Essays presented in Honor of Karel Lambert}. Dordrecht,
  Boston, London: Kluwer.

\end{thebibliography}
\end{document}